\documentclass[a4paper,twoside]{article}
\usepackage{amssymb,amsmath,amscd,amsthm}

\newtheorem{theorem}{Theorem}[section]
\newtheorem{lemma}[theorem]{Lemma}
\newtheorem{corollary}[theorem]{Corollary}
\newtheorem{proposition}[theorem]{Proposition}

\theoremstyle{definition}
\newtheorem{definition}[theorem]{Definition}
\newtheorem{example}[theorem]{Example}

\theoremstyle{remark}
\newtheorem{remark}[theorem]{Remark}

\newtheorem*{acknowledgments}{Acknowledgments}

\newcommand{\pro}[2]{\langle #1, #2 \rangle}

\def\lolra{\Longleftrightarrow}
\def\lra{{\,\Leftrightarrow\,}}
\def\ra{{\,\Rightarrow\,}}

\def\N{{\mathbb N}}
\def\Z{{\mathbb Z}}
\def\Q{{\mathbb Q}}
\def\R{{\mathbb R}}

\def\F{{\mathcal{F}}}
\def\V{{{\mathcal{V}}}}

\def\Vs{{\tilde{V}}}

\def\MR{{M_{\R}}}
\def\NR{{N_{\R}}}

\def\conv{{\rm conv}}

\def\Hom{{\rm Hom}}
\def\Mat{{\rm Mat}}
\def\GL{{\rm GL}}
\def\id{{\rm id}}

\def\Herm{{\rm Herm}}

\def\rand{\partial}
\def\randp{{\rand P}}

\numberwithin{equation}{section}

\newcommand{\abs}[1]{\lvert#1\rvert}

\begin{document}

\title{Classification of pseudo-symmetric simplicial reflexive polytopes}

\author{{\sc Benjamin Nill\vspace{0.5ex}} \\
\small  {\em Research Group Lattice Polytopes, FU Berlin}   \\
\small  {\em Arnimallee 3, 14195 Berlin}  \\
\small  {\em e-mail: nill@math.fu-berlin.de} \\
 }

\date{{\normalsize November 2005}}

\maketitle

\begin{abstract}
Gorenstein toric Fano varieties correspond to so-called reflexive polytopes. 
If such a polytope contains a centrally symmetric pair of facets, we call the polytope, respectively 
the toric variety, pseudo-symmetric. Here we present a complete classification of 
pseudo-symmetric simplicial reflexive polytopes. This is a generalization of a result of Ewald on 
pseudo-symmetric nonsingular toric Fano varieties and recent work of Wirth. 
As applications we determine the maximal number of vertices, facets and lattice points, and show that the vertices can be chosen 
to have coordinates $-1,0,1$.
\end{abstract}

\section*{Introduction}

Isomorphism classes of nonsingular toric Fano varieties correspond to 
unimodular isomorphism classes of so-called {\em smooth Fano polytopes}; these are lattice polytopes, where 
the origin $0$ of the lattice is contained in their interiors, and where the vertices of any facet 
form a basis of the lattice. 
Smooth Fano polytopes have been classified up to dimension four, see \cite{WW82,Bat82,Bat99,Sat00}. 
However in higher dimensions classification results require more symmetries of the polytope. 

Let us call a toric Fano variety 
{\em centrally symmetric}, if the associated Fano polytope is centrally symmetric with respect to $0$. 
Moreover, let us denote by a {\em del Pezzo polytope} $V_d$ the $d$-dimensional 
centrally symmetric smooth Fano polytope with $2d+2$ vertices 
$\pm e_1, \ldots, \pm e_d, \pm (e_1 + \cdots + e_d)$, where $e_1, \ldots, e_d$ is a basis of the lattice 
of even rank $d$. We call the corresponding $d$-dimensional nonsingular toric Fano variety a {\em Voskresenskij-Klyachko variety} 
(previously called del Pezzo variety). In \cite{VK85} Voskresenskij and Klyachko showed that 
{\em any centrally symmetric nonsingular toric Fano variety 
is a product of projective lines and Voskresenskij-Klyachko varieties.}

In \cite{Ewa88} Ewald gave a 
generalization of this result by assuming that the polytope is only pseudo-symmetric. 
We say that a polytope $P$ with interior point $0$, 
respectively the associated variety, is {\em pseudo-symmetric}, if 
$P$ has a centrally symmetric pair of facets. We can define a 
pseudo-symmetric smooth Fano polytope $\Vs_d$ called {\em pseudo-del Pezzo polytope} 
as the convex hull of the $2d + 1$ vertices $\pm e_1, \ldots, \pm e_d, -e_1 - \cdots - e_d$ 
in above notation. 
We call the associated variety an {\em Ewald variety} (previously called pseudo-del Pezzo variety). 
Now in \cite{Ewa88} Ewald showed that 
{\em any pseudo-symmetric nonsingular toric Fano variety 
is a product of projective lines, Voskresenskij-Kly\-achko varieties and Ewald varieties.} 
In convex-geometric language this means that any pseudo-symmetric smooth Fano polytope $P$ {\em splits} into copies of $[-1,1]$, del Pezzo polytopes 
and pseudo-del Pezzo polytopes. By a result of Casagrande in \cite{Cas03} this holds also for any $d$-dimensional smooth Fano polytope $P$ 
having $d$ linearly independent vertices $v_1, \ldots, v_d$ such that $-v_1, \ldots, -v_d$ are also vertices in $P$.

\smallskip

In the context of mirror symmetry Batyrev introduced in \cite{Bat94} the notion of a reflexive polytope 
that is weaker than that of a smooth Fano polytope:  A fully-dimensional lattice polytope containing 
the origin $0$ in its interior is called {\em reflexive},  
if for any facet $F$ the unique vector in the dual vector space that evaluates $-1$ on $F$ 
is a lattice point. This implies that the dual polytope is also a reflexive polytope. 
Isomorphism classes of reflexive polytopes correspond 
to toric Fano varieties having at most Gorenstein singularities, i.e., projective toric varieties whose anticanonical divisor 
is an ample Cartier divisor. 

Reflexive polytopes are interesting classes of lattice polytopes 
regarded from several different aspects and 
were classified up to dimension four using computer algorithms, see \cite{KS98,KS00, KS04}. 
In this paper we present as a rare higher-dimensional classification result a 
generalization of Ewald's theorem by regarding not only smooth Fano polytopes 
but simplicial reflexive polytopes, i.e., reflexive polytopes where any facet is a simplex. 
This is indeed a significant extension, since 
for instance in dimension four there are $124$ isomorphism classes of smooth Fano polytopes \cite{Bat99, Sat00} 
compared to $5450$ isomorphism classes of simplicial reflexive polytopes in the database \cite{KS05}. 
In convex-geometric language our main result reads as follows, 
for this recall that a crosspolytope is the combinatorial dual of a cube:

\begin{theorem}
Any pseudo-symmetric simplicial reflexive polytope splits up to unimodular isomorphisms uniquely into 
a centrally symmetric reflexive crosspolytope, del Pezzo polytopes, and pseudo-del Pezzo polytopes. 
The isomorphism class of any centrally symmetric reflexive crosspolytope in fixed dimension 
can be determined by a finite number of suitable matrix normal forms.
\label{intro}
\end{theorem}

A more precise version of this theorem will be given in section two, divided into 
Theorem \ref{theo1} and Theorem \ref{generalewald}. 
In the case of a smooth Fano polytope our main result immediately yields the theorem of Ewald in 
\cite{Ewa88}. The second part of the theorem was already formulated and 
proven by Wirth, a student of Ewald, 
in \cite{Wir97} using the Hermite normal form theorem that we also rely on here. 
However our proof is independent of the results in Ewald and Wirth. While their proofs depend on 
explicit determinant calculations, we focus on general discrete-geometric properties 
of simplicial reflexive polytopes and 
rather deal with dual bases and the dual reflexive polytope. For this we apply observations 
on reflexive polytopes from \cite{Nil05}.

\smallskip

Some applications:

\begin{itemize}
\item We explicitly carry out the classification of all pseudo-symmetric simplicial 
reflexive polytopes up to dimension six.
\item In any dimension we show that there is up to isomorphism only one pseudo-symmetric simplicial 
reflexive polytope with the maximal number of vertices. The same statements holds for 
the maximal number of lattice points.
\item The dual polytope of a $d$-dimensional pseudo-symmetric simplicial reflexive polytope $P$ 
has at most $6^{d/2}$ vertices, where the equality case is only achieved, if $P$ splits into $d/2$ copies of $V_2$.
\end{itemize}

The last point is especially interesting, since we get a precise confirmation of the 
general conjecture \cite[Conjecture 5.2]{Nil05} 
on the maximal number of vertices of a reflexive polytope.

Finally turning to a more well-known conjecture, in \cite{Ewa88} Ewald 
conjectured that one can always embed any $d$-dimensional smooth Fano polytope 
in $[-1,1]^d$, as observed for pseudo-symmetric smooth Fano polytopes. This conjecture is no longer valid 
for general simplicial reflexive polytopes, however we show that it still holds in the pseudo-symmetric 
case, and moreover that their dual polytopes can be embedded into  $\lfloor\frac{d}{2}\rfloor [-1,1]^d$.

\smallskip

This article is organized in the following way: 
In the first section the notation is fixed and basic notions are recalled. 
In the second section a refined version of the main result is formulated, and proven in the third section. 
The last section contains applications.

\smallskip

{\em Remark:} On advice of Batyrev we have avoided the term 'del Pezzo variety', previously used in \cite{VK85,Ewa88,Ewa96,Cas03}, 
since there is by now an established theory of del Pezzo varieties in the sense of Fujita.

\section{Basic notions and preliminary results}

\subsection*{The notation}

Here we set up the setting of this paper:

\begin{itemize}
\item $M \cong \Z^d$ and $N := \Hom_\Z(M,\Z)$ are {\em dual lattices} 
with the pairing $\pro{\cdot}{\cdot}$. We set $\MR := M \otimes_\Z \R \cong \R^d$ and 
$\NR := N \otimes_\Z \R \cong \R^d$.
\item For a subset $S \subseteq \MR$ we denote by $\conv(S)$ the convex hull of $S$, and by $\dim(S)$ its dimension.
\item A {\em lattice polytope} in $\MR$ is the convex hull of lattice points in $M$. 
Always let $P$ be a $d$-dimensional lattice polytope in $\MR$ 
that contains the origin $0$ in its interior and whose vertices are primitive lattice points. 
Such a $P$ is called {\em Fano polytope}.
\item The set of vertices of $P$ is denoted by $\V(P)$, the set of facets by $\F(P)$. 
The boundary of $P$ is refered to as $\randp$.
For any facet $F \in \F(P)$ there is a unique {\em inner normal} $\eta_F \in \NR$ defined by 
$\pro{\eta_F}{F} = -1$.
\item Two lattice polytopes are regarded as {\em isomorphic}, if they are isomorphic 
under some unimodular transformation, 
i.e., if there is a lattice automorphism of $M$ that maps the vertex sets mutually 
onto each other. 
\item $P$ spans a fan over its faces, denoted by $\Sigma_P$. The associated {\em toric Fano variety} 
is denoted by $X(M,\Sigma_P)$, briefly $X(\Sigma_P)$. This sets up a 
one-to-one-correspondence between isomorphism classes of Fano polytopes and 
isomorphism classes of toric Fano varieties, see \cite{Bat94} or \cite{Nil05}.
\end{itemize}

\subsection*{Reflexive polytopes}

The {\em dual polytope} of $P$ is defined as 
$$P^* := \{x \in \NR \,:\, \pro{x}{y} \geq -1 \;\forall\, y \in P\}$$
and has as vertices precisely the inner normals of the facets of $P$. 

\begin{definition}
\begin{itemize}
\item[]
\item $P$ is called {\em smooth Fano polytope}, if the vertices of any facet form a lattice basis of $M$. 
Equivalently, $X(\Sigma_P)$ is a nonsingular toric Fano variety.
\item $P$ is called {\em reflexive polytope}, if $P^*$ is a lattice polytope. 
Equivalently, $X(\Sigma_P)$ is a toric Fano variety with at most Gorenstein singularities.
\end{itemize}
\end{definition}

There is the following duality of reflexive polytopes, see \cite{Bat94}:

\begin{center}
$P \subseteq \MR$ is reflexive $\,\lolra\,$ $P^* \subseteq \NR$ is reflexive.
\end{center}

Furthermore we need a simple but fundamental property \cite[Prop. 4.1]{Nil05}:

\begin{lemma}
Let $P \subseteq \MR$ be a reflexive polytope, and $v, w \in \randp \cap M$.

If $v + w \not= 0$ and there is no facet containing both $v$ and $w$, then 
$v + w \in \randp \cap M$.
\label{prim}
\end{lemma}

Another important observation is the following result, that is included in 
\cite[Lemma 5.5]{Nil05}:

\begin{lemma}
Let $P \subseteq \MR$ be a reflexive polytope, $F \in \F(P)$, and $m \in \randp \cap M$.

If $\pro{\eta_F}{m} = 0$, then $m$ is contained in a facet intersecting $F$ in a codimension two face.
\label{adjacent}
\end{lemma}

\subsection*{Pseudo-symmetry and del Pezzo polytopes}

\begin{definition}
\begin{itemize}
\item[]
\item $P$ is called {\em centrally symmetric}, if $-P = P$.
\item $P$ is called {\em pseudo-symmetric}, if there is some $F \in \F(P)$ with $-F \in \F(P)$.
\end{itemize}
\label{delpezzos}
\end{definition}

For us the following smooth Fano polytopes will be especially important:

\begin{definition}
Let $e_1, \ldots, e_d$ be a lattice basis of $M$. Let $d$ be even.
\begin{itemize}
\item $V_d := \conv(\pm e_1, \ldots, \pm e_d, \pm (e_1 + \cdots + e_d))$ 
is called a {\em del Pezzo polytope}.
\item $\Vs_d := \conv(\pm e_1, \ldots, \pm e_d, -e_1 - \cdots - e_d)$ 
is called a {\em pseudo-del Pezzo polytope}.
\end{itemize}
\end{definition}

So del Pezzo polytopes are centrally symmetric, while pseudo-del Pezzo polytopes 
are only pseudo-symmetric.

\subsection*{Hermite matrix normal forms}

Let $n \in \N_{\geq 1}$.

\begin{definition}
\begin{itemize}
\item[]
\item $\GL_n(\Z)$ is the set of $n \times n$-matrices with integer coefficients and 
determinant $\pm 1$.
\item For an arbitrary set of integers $R$ we 
let $\Mat_{m \times n}(R)$ be the set of $m \times n$-matrices with entries in $R$. 
We abbreviate $\Mat_n(R) := \Mat_{n \times n}(R)$. 
\item For $n, \lambda \in \N_{\geq 1}$ we denote by $\Herm(n,\lambda)$ the finite set 
of lower triangular matrices $H \in \Mat_n(\N)$ 
with determinant $\lambda$ satisfying $h_{i,j} < h_{j,j}$ for all $j = 1, \ldots, n-1$ and $i > j$.
\end{itemize}
\end{definition}

The famous theorem of Hermite is the following (e.g., \cite[pp. 15-18]{New72}): 

\begin{theorem}
For any $L \in \Mat_n(\Z)$ with determinant $\lambda \not= 0$ there exist matrices $U \in \GL_n(\Z)$ and 
$H \in \Herm(n,\lambda)$ such that $U L = H$.
\label{hermitethm}
\end{theorem}

\section{The main theorems}

In this section Theorem \ref{intro} will be formulated more precisely. 
It is split into Theorem \ref{theo1} and Theorem \ref{generalewald}, the first one dealing 
with the case of the minimal number of vertices.

\subsection*{Classification of centrally symmetric reflexive crosspolytopes}

Centrally symmetric reflexive crosspolytopes have been classified in \cite[Satz 3.3]{Wir97}. Here we state Wirth's 
result in a somewhat strengthened form.

To simplify notation we define:

\begin{definition}
\begin{itemize}
\item[]
\item A matrix $A \in \Mat_d(\N)$ is called {\em Wirth matrix}, if 
$$A = \begin{pmatrix}2 \id_{f} & 0\\ C & \id_{d-f}\end{pmatrix},$$
where $f \in \{0, \ldots, d-1\}$ and $C \in \Mat_{(d-f) \times f}(\{0,1\})$ such that any column of $C$ 
has an odd number of $1$'s. Here $\id_k$ is the $k \times k$-identity matrix.

\item A Wirth matrix $A$ is called {\em 1-minimal Wirth matrix}, if any row of $C$ contains some $1$. 
Obviously we get from a Wirth matrix a 1-minimal one called its {\em reduction} 
by deleting rows containing only one $1$ (and the corresponding columns).

\item Two matrices in $\Mat_d(\N)$ are regarded as {\em equivalent}, if they differ only up to 
permutation of columns and left-multiplication by a matrix in $\GL_d(\Z)$.

\item Let $P_1, P_2$ be non-zero Fano polytopes, with respect to lattices $M_1$, $M_2$, such that 
$\dim(P_1) + \dim(P_2) = d$. We say $P$ {\em splits} into factors $P_1$ and $P_2$, if 
$P \cong \conv(P_1 \oplus \{0\}, \{0\} \oplus P_2)$ for $M \cong M_1 \oplus M_2$. 
Equivalently, $X(\Sigma_P) \cong X(\Sigma_{P_1}) \times X(\Sigma_{P_2})$; or dually, $P^* \cong P_1^* \times P_2^*$.

\item We say $P$ is {\em 1-irreducible}, if $P$ does not split into $[-1,1]$ and some $P_2$.
\item We denote by a {\em cs-crosspolytope} a centrally symmetric crosspolytope, i.e., 
the convex hull of a simplex (not containing $0$) and its negative. Equivalently, 
a $d$-dimensional cs-crosspolytope is a $d$-dimensional pseudo-symmetric simplicial polytope with $2 d$ vertices.
\end{itemize}
\label{Wirth}
\end{definition}

Obviously any matrix in $\Mat_d(\Z)$ with non-zero determinant defines a lattice cs-crosspolytope 
by taking the convex hull of its columns and their negatives. 

Thereby reflexive cs-crosspolytopes can be classified:

\begin{theorem}[Wirth, N.]
\begin{itemize}
\item[]
\item There is a one-to-one correspondence between equivalence classes of Wirth matrices and 
isomorphism classes of reflexive cs-crosspolytopes. For the inverse map take any facet of a reflexive cs-crosspolytope, 
then there is a lattice basis 
such that the coordinates of the vertices of this facet are the columns of a corresponding Wirth matrix. 
\item Hereby equivalence classes of {\em 1-minimal} Wirth matrices 
correspond to isomorphism classes of {\em 1-irreducible} reflexive cs-crosspolytopes. 
Moreover any reflexive cs-crosspolytope $P$ splits up to isomorphism uniquely into a 1-irreducible 
reflexive cs-cross\-polytope $P'$ and $r$ copies of $[-1,1]$. Here an associated 
1-minimal Wirth matrix of $P'$ is just 
the reduction of a Wirth matrix $A$ associated to $P$, and $r$ equals the number of rows of $A$ 
containing only one $1$.
\end{itemize}
\label{theo1}
\end{theorem}

To sum up, the isomorphism class of any reflexive cs-crosspolytope $P$ is given by an up to equivalence 
uniquely determined 1-minimal Wirth matrix $A'$ (encoding the {\em singular} factor $P'$ of $P$) 
and a unique natural number $r$ (corresponding to the {\em nonsingular} factor of $P$). 
The determinant of $A'$ is just the index in $M$ of the lattice 
spanned by the vertices of $P$.

\begin{example}
As an illustration we give a list of all isomorphism classes of 1-irreducible reflexive cs-crosspolytopes by 
listing the equivalence classes of their associated 1-minimal Wirth matrices for $d \leq 6$ (this list is for $d \leq 5$ implicitly already 
contained in \cite[3.7,3.8]{Wir97}):

\begin{itemize}
\item[] $d=2$:
$\begin{pmatrix}2 & 0\\1 & 1\end{pmatrix}$
\item[] $d=3$:
$\begin{pmatrix}2 & 0 & 0\\0 & 2 & 0\\1 & 1 & 1\end{pmatrix}$
\item[] $d=4$:
$\begin{pmatrix}2 & 0 & 0 & 0\\1 & 1 & 0 & 0\\1 & 0 & 1 & 0\\1 & 0 & 0 & 1\end{pmatrix}$, 
$\begin{pmatrix}2 & 0 & 0 & 0\\0 & 2 & 0 & 0\\1 & 0 & 1 & 0\\0 & 1 & 0 & 1\end{pmatrix}$, 
$\begin{pmatrix}2 & 0 & 0 & 0\\0 & 2 & 0 & 0\\0 & 0 & 2 & 0\\1 & 1 & 1 & 1\end{pmatrix}$
\item[] $d=5$:
$\begin{pmatrix}2 & 0 & 0 & 0 & 0\\0 & 2 & 0 & 0 & 0\\1 & 1 & 1 & 0 & 0\\ 1 & 1 & 0 & 1 & 0\\ 
1 & 1 & 0 & 0 & 1\end{pmatrix}$,
$\begin{pmatrix}2 & 0 & 0 & 0 & 0\\0 & 2 & 0 & 0 & 0\\0 & 0 & 2 & 0 & 0\\1 & 0 & 0 & 1 & 0\\0 & 1 & 1 & 0 & 1\end{pmatrix}$, 
$\begin{pmatrix}2 & 0 & 0 & 0 & 0\\0 & 2 & 0 & 0 & 0\\0 & 0 & 2 & 0 & 0\\0 & 0 & 0 & 2 & 0\\1 & 1 & 1 & 1 & 1
\end{pmatrix}$
\item[] $d=6$: 
$\left(\begin{smallmatrix}2 & 0 & 0 & 0 & 0 & 0\\1 & 1 & 0 & 0 & 0 & 0\\1 & 0 & 1 & 0 & 0 & 0\\
1 & 0 & 0 & 1 & 0 & 0\\1 & 0 & 0 & 0 & 1 & 0\\1 & 0 & 0 & 0 & 0 & 1\end{smallmatrix}\right)$,
$\left(\begin{smallmatrix}2 & 0 & 0 & 0 & 0 & 0\\0 & 2 & 0 & 0 & 0 & 0\\1 & 0 & 1 & 0 & 0 & 0\\
1 & 0 & 0 & 1 & 0 & 0\\1 & 0 & 0 & 0 & 1 & 0\\0 & 1 & 0 & 0 & 0 & 1\end{smallmatrix}\right)$
$\left(\begin{smallmatrix}2 & 0 & 0 & 0 & 0 & 0\\0 & 2 & 0 & 0 & 0 & 0\\1 & 1 & 1 & 0 & 0 & 0\\
1 & 1 & 0 & 1 & 0 & 0\\1 & 0 & 0 & 0 & 1 & 0\\0 & 1 & 0 & 0 & 0 & 1\end{smallmatrix}\right)$
$\left(\begin{smallmatrix}2 & 0 & 0 & 0 & 0 & 0\\0 & 2 & 0 & 0 & 0 & 0\\0 & 0 & 2 & 0 & 0 & 0\\
1 & 0 & 0 & 1 & 0 & 0\\0 & 1 & 0 & 0 & 1 & 0\\0 & 0 & 1 & 0 & 0 & 1\end{smallmatrix}\right)$,\\
$\left(\begin{smallmatrix}2 & 0 & 0 & 0 & 0 & 0\\0 & 2 & 0 & 0 & 0 & 0\\0 & 0 & 2 & 0 & 0 & 0\\
1 & 1 & 1 & 1 & 0 & 0\\1 & 0 & 0 & 0 & 1 & 0\\1 & 0 & 0 & 0 & 0 & 1\end{smallmatrix}\right)$,
$\left(\begin{smallmatrix}2 & 0 & 0 & 0 & 0 & 0\\0 & 2 & 0 & 0 & 0 & 0\\0 & 0 & 2 & 0 & 0 & 0\\
1 & 1 & 1 & 1 & 0 & 0\\1 & 1 & 0 & 0 & 1 & 0\\1 & 1 & 0 & 0 & 0 & 1\end{smallmatrix}\right)$, 
$\left(\begin{smallmatrix}2 & 0 & 0 & 0 & 0 & 0\\0 & 2 & 0 & 0 & 0 & 0\\0 & 0 & 2 & 0 & 0 & 0\\
0 & 0 & 0 & 2 & 0 & 0\\1 & 1 & 1 & 0 & 1 & 0\\0 & 0 & 0 & 1 & 0 & 1\end{smallmatrix}\right)$, 
$\left(\begin{smallmatrix}2 & 0 & 0 & 0 & 0 & 0\\0 & 2 & 0 & 0 & 0 & 0\\0 & 0 & 2 & 0 & 0 & 0\\
0 & 0 & 0 & 2 & 0 & 0\\1 & 1 & 0 & 0 & 1 & 0\\0 & 0 & 1 & 1 & 0 & 1\end{smallmatrix}\right)$, 
$\left(\begin{smallmatrix}2 & 0 & 0 & 0 & 0 & 0\\0 & 2 & 0 & 0 & 0 & 0\\0 & 0 & 2 & 0 & 0 & 0\\
0 & 0 & 0 & 2 & 0 & 0\\0 & 0 & 0 & 0 & 2 & 0\\1 & 1 & 1 & 1 & 1 & 1\end{smallmatrix}\right)$\\
\end{itemize}
\smallskip

Note that for $d=4,5,6$ there are $1,1,4$ of these Wirth matrices whose corresponding 
reflexive cs-crosspolytope splits into smaller factors.
\label{list}
\end{example}

\begin{remark}
Two matrices $A_1$ and $A_2$ in $\Mat_d(\Z)$ with the same determinant are equivalent, if for some permutation $\pi$ 
all entries of $A_1^\pi A_2^{-1}$ are integers, where $A_1^\pi$ is the matrix of $\pi$-permuted 
columns of $A_1$. For Wirth matrices this can be easily checked due to their simple structure, see \cite[Satz 3.9]{Wir97}.
\end{remark}

\subsection*{Classification of pseudo-symmetric simplicial reflexive polytopes} 

\begin{theorem}
Let $P \subseteq \MR$ be a pseudo-symmetric simplicial reflexive polytope. Then 
$P$ splits up to isomorphism uniquely into a 1-irreducible reflexive cs-crosspolytope 
$P'$, $r$ copies of $[-1,1]$, del Pezzo polytopes, and pseudo-del Pezzo polytopes.
\label{generalewald}
\end{theorem}

So $P$ splits uniquely into $P'$ (the {\em singular} factor) and a smooth Fano polytope (the {\em nonsingular} factor). 
We recover the original result of Ewald in \cite{Ewa88} under milder assumptions:

\begin{corollary}
Let $P \subseteq \MR$ be a pseudo-symmetric simplicial reflexive polytope where the vertices span 
the lattice $M$.

Then the corresponding toric variety $X(\Sigma_P)$ is just 
a product of projective lines, Voskresenskij-Klyachko varieties, and Ewald varieties. 
In particular $X(\Sigma_P)$ is nonsingular and $P$ a smooth Fano polytope.
\label{ewald}
\end{corollary}

\begin{proof}

Indeed, since by assumption the determinant of the 1-minimal Wirth matrix associated to $P'$ equals one, 
$P'$ has to be zero.

\end{proof}

Now we can easily calculate that there are $1,3,3,8,8,18$ isomorphism classes of $d=1,2,3,4,5,6$-dimensional 
pseudo-symmetric smooth Fano polytopes, therefore Theorem \ref{generalewald} together with the list in Example \ref{list} yields:

\begin{corollary}
For $d=2,3,4,5,6$ there are exactly $4,5,15,20,50$ isomorphism classes of $d$-dimensional 
pseudo-symmetric simplicial reflexive polytopes.
\end{corollary}

A rigorous proof of the previous result was in the centrally symmetric case 
up to now only known for $d \leq 3$, see \cite{Wag95}. For $d=4$ the centrally symmetric case of up to $10$ 
vertices was dealt with by rather complicated and long calculations in 
\cite[Satz 5.11]{Wir97}. 

\smallskip

Applying \cite[Cor. 1.16]{Oda88} (and, e.g., \cite[Prop. 1.9]{Nil05}) to \ref{generalewald} yields also:

\begin{corollary}
Any $d$-dimensional pseudo-symmetric $\Q$-factorial Gorenstein toric Fano variety 
is the projection for the quotient of a product of 
projective lines, Voskresenskij-Klyachko varieties, and Ewald varieties 
with respect to the action of a finite group isomorphic to 
$(\Z/2\Z)^f$ for $f \leq d-1$.
\label{csquot}
\end{corollary}

The combinatorial statement sounds rather surprising:

\begin{corollary}
Any pseudo-symmetric simplicial reflexive polytope is combinatorially isomorphic to 
a pseudo-symmetric smooth Fano polytope.
\label{combcorosimpliz}
\end{corollary}

\begin{remark}
One should not be misled by this result: Without the symmetry assumption the combinatorics 
of simplicial reflexive polytopes can be much more complicated than the one of smooth Fano polytopes. 
For instance, according to the database \cite{KS05} 
the columns of the following matrix form the vertices of the only four-dimensional reflexive polytope with $7$ vertices 
and $14$ facets; it is simplicial but not a smooth Fano polytope: $\left(\begin{smallmatrix}
    1 &   0 &   0  &  0  &-2  &  \text{ }0 &   1\\
    0 &   1 &   0  &  0  &  \text{ }0  &-2 &   1\\
    0 &   0 &   1  &  0  &-1  &  \text{ }0 &   1\\
    0 &   0 &   0  &  1  &  \text{ }0  &-1 &   1
\end{smallmatrix}\right)$
\end{remark}

\section{Proof of the main theorems}

In this section Theorems \ref{theo1} and \ref{generalewald} will be proven. 
We will deal directly with the general case of a pseudo-symmetric simplicial reflexive polytope, 
this includes therefore an independent proof of the results in \cite{Wir97} on reflexive cs-crosspolytopes.

\subsection*{The key-lemma}

As a preparation we need the following important fact (pay attention, in what follows 
$e_1, \ldots, e_d$ is in general {\em not} a lattice basis):

\begin{lemma}
Let $P$ be a simplicial reflexive polytope with facets $F$, $-F$. 

Let $\V(F) = \{e_1, \ldots, e_d\}$, and $e_1^*, \ldots, e_d^*$ be the dual $\R$-basis of $\NR$. 
For $i = 1, \ldots, d$ we denote by $F_i$ the unique facet of $P$ 
such that $F_i \cap F = \conv(e_j \,:\, j \not=i)$. We set $u := \eta_F \in \V(P^*)$.  

Let $v \in \V(P) \cap u^\perp$. We write $v = \sum_{i=1}^d q_i e_i$ as a rational linear combination. 
Then we have for $i = 1, \ldots, d$:
\[q_i < 0 \,\lolra\, q_i = -1 \,\lolra\, v \in F_i.\]
In this case $e_i^* = \eta_{F_i} - u \in P^* \cap N$.\\

Moreover there are $I, J \subseteq \{1, \ldots, d\}$ with $I \cap J = \emptyset$ and $\abs{I} = \abs{J}$ 
such that
\[v = \sum_{j \in J} e_j - \sum_{i \in I} e_i.\]
\label{fund2}
\end{lemma}

\begin{proof}

In the proof of \cite[Lemma 5.5]{Nil05} it was already readily observed that 
$\eta_{F_i} = u + \alpha_i e_i^*$ for $\alpha_i > 0$ and $i = 1, \ldots, d$. 
Hence we get $q_i = \pro{e_i^*}{v} < 0 \,\lra\, v \in F_i$. 

Now let $q_i = \pro{e_i^*}{v} < 0$, $v \in F_i$. 
Therefore $F_i \cap (-F) = \emptyset$, so by duality we can apply Lemma \ref{prim} to the reflexive polytope 
$P^* \subseteq \NR$ and the inner normals of $F_i$ and $-F$. This yields $\eta_{F_i} - u = \alpha_i e_i^* 
\in P^* \cap N$. Hence $-1 \leq \pro{\alpha_i e_i^*}{-e_i} = - \alpha_i \in \Z$, so $\alpha_i = 1$, 
and $q_i = \pro{e_i^*}{v} = \pro{\eta_{F_i} - u}{v} = -1$. 

Obviously $v = \sum_{j=1}^d (-q_j) (-e_j)$, 
so applying the previous result to $-F$ yields $q_j > 0 \,\lra\, q_j = 1$, hence the remaining statement 
due to $\pro{u}{v} = 0$.

\end{proof}

\subsection*{The structure theorem}

\begin{proposition}
Let $P$ be a simplicial reflexive polytope with facets $F, -F$.

Let $\V(F) = \{e_1, \ldots, e_d\}$. For $i = 1, \ldots, d$ we let $v^i$ denote 
the unique vertex of $P$ contained in the unique facet that intersects $F$ in 
$\conv(e_j \,:\, j \not= i)$.

\begin{itemize}
\item There exists a lattice basis $m_1, \ldots, m_d$ of $M$ such that in this basis the coefficient vectors of $e_1, \ldots, e_d$ 
are the columns of a Wirth matrix $A$.
\item Any vertex of $P$ is in $\{\pm e_1, \ldots, \pm e_d, v^1, \ldots, v^d\}$. 
\item There exist pairwise disjoint subsets $I_1, \ldots, I_l \subseteq \{1, \ldots, d\}$ 
and pairwise disjoint subsets $J_1, \ldots, J_l \subseteq \{1, \ldots, d\}$ 
with $I_k \cap J_k = \emptyset$ and $\abs{I_k} = \abs{J_k}$ 
for all $k = 1, \ldots, l$ such that for $i \in \{1, \ldots, d\}$ we have 
\[v^i \in \eta_F^\perp \;\lolra\; i \in \bigcup_{k = 1}^l I_k,\]
and moreover, if $i \in I_k$ for some $k \in \{1, \ldots, l\}$, then
\[v^i = \sum_{j \in J_k} e_j - \sum_{i' \in I_k} e_{i'}.\]
\item For $i \in I_1 \cup \cdots I_l \cup J_1 \cup \cdots \cup J_l$ the $i$th row of $A$ 
is of the form \newline $(0, \ldots, 0, 1, 0, \ldots, 0)$ for $1$ at the $i$th position. 
\item If for $k,k' \in \{1, \ldots, l\}$ the sets $I_k$ and $J_{k'}$ intersect, then $k \not= k'$, 
$I_k = J_{k'}$, and $J_k = I_{k'}$.
\end{itemize}
\label{simplizchar}
\end{proposition}

\begin{proof}

Let $u$, $F_i$ defined as in Lemma \ref{fund2}. Consider the following steps 
for the construction of $A$ and $m_1, \ldots, m_d$:

\begin{enumerate}
\item Let $i \in \{1, \ldots, d\}$. 

If $v^i \in u^\perp$, then by Lemma \ref{fund2} $e^*_i = \eta_{F_i} - u \in N$.

If $v^i \not\in u^\perp$, then obviously $v^i = - e_i$, so $e^*_i = \frac{\eta_{F_i} - u}{2} \in \frac{1}{2} N$.

Hence in any case $e^*_1, \ldots, e^*_d \in \frac{1}{2} N$. In particular this yields 
$$2 M \subseteq \Z e_1 + \cdots + \Z e_d \subseteq M.$$

\item We define for an arbitrary lattice basis of $M$ the matrix $L \in \Mat_d(\Z) \cap \GL_d(\Q)$, 
where the columns are the coefficient vectors of $e_1, \ldots, e_d$ in this basis. By Theorem \ref{hermitethm} 
there exists $U \in \GL_d(\Z)$ such that $A := U L \in \Mat_d(\N)$ is a lower triangular matrix 
with $A_{i,j} \in \{0, \ldots, A_{j,j} - 1\}$ for $i > j$. Hence there is a lattice basis $m_1, \ldots, m_d$ of $M$ 
such that the columns of $A$ are the coefficient vectors of $e_1, \ldots, e_d$ in the lattice basis $m_1, \ldots, m_d$.

\item The first point yields that $2 m_i$ is contained for $i = 1, \ldots, d$ in the column space of $A$, in particular 
any diagonal element of $A$ equals $1$ or $2$. 

\item If $A_{i,j} = 1$ for $i > j$, then necessarily $A_{j,j} = 2$ and $A_{i,i} = 1$ (again one has to use that 
$2 m_j$ is a $\Z$-linear combination of $e_j, \ldots, e_d$.)

\item Using the previous point we assume, by possibly permutating the columns and the rows of $A$, 
that $A$ has a blockmatrix structure as in Definition \ref{Wirth}. Since any vertex of a reflexive polytope is primitive, 
obviously $f \not= d$. It remains to show that any column has an odd number of $1$'s: 
By the previous point let $e_j = 2 m_j + \sum_{k=1}^s e_{i_k}$, where $i_k > j$. 
We get $2 \pro{\eta_F}{m_j} =  \pro{\eta_F}{e_j - \sum_{k=1}^s e_{i_k}} = -1 + s$, so $s$ has to be odd. 
Hence $A$ is a Wirth matrix.

\end{enumerate}

Now it is easy to see that 
$e_i^* \in N$ if and only if the $i$th row of $A$ is of the form $(0, \ldots, 0, 1, 0, \ldots, 0)$. 
By this observation, Lemma \ref{fund2} (applied to $F$ and $-F$) 
and the first point in the proof we see that it only remains to show the last statement in the proposition:

So let $k,k' \in \{1, \ldots, l\}$ with $I_k \cap J_{k'} \not= \emptyset$. By construction 
$I_k \cap J_k = \emptyset$, hence $k \not= k'$. 

Assume $I_{k'} \not\subseteq J_k$. 
Then there exists $i \in I_{k'}$, $i \not\in J_k$. Let $j \in I_k \cap J_{k'}$, in particular $j \not= i$. 
Now we define for the dual $\R$-basis $e_1^*, \ldots, e_d^*$ of $\NR$ the vector 
$w := e_i^* - \sum_{s \not= i,j} e_s^*$. By construction it is easy to check that 
$w$ is an inner normal of a face of $P$ containing as vertices 
$e_s$ for $s = 1, \ldots, d$ with 
$s \not= i,j$, as well as $-e_i$ and $v^i$ and $v^j$. This is a contradiction to 
$P$ being simplicial.

Hence $I_{k'} \subseteq J_k$. In particular $I_{k'} \cap J_k \not= \emptyset$, so also 
$I_k \subseteq J_{k'}$. Since therefore $\abs{I_{k'}} \leq \abs{J_k} = \abs{I_k} \leq \abs{J_{k'}} = \abs{I_{k'}}$, 
we have $I_{k'} = J_k$ and $I_k = J_{k'}$.

\end{proof}

\subsection*{Proof of Theorem \ref{theo1}}

By the first point of Proposition \ref{simplizchar} any reflexive cs-crosspolytope is defined by 
a Wirth matrix. On the other hand it was shown in \cite{Wir97} by a simple calculation 
that any Wirth matrix defines a reflexive cs-crosspolytope. Furthermore if two Wirth matrices define 
isomorphic reflexive cs-crosspolytopes, then to see that the Wirth matrices have to be equivalent, it 
is enough to show that any two facets of a reflexive cs-crosspolytope are isomorphic as lattice polytopes. 

For this let us assume that $P$ is a reflexive cs-crosspolytope defined by a Wirth 
matrix $A$ with columns $e_1, \ldots, e_d$. Now we regard a matrix where the columns are given as 
$\epsilon_1 e_1, \ldots, \epsilon_d e_d$ for $\epsilon_1, \ldots, \epsilon_d \in \{-1,1\}$. 
We perform elementary row operations on this matrix by first 
multiplying any row that has a negative number on the diagonal by $-1$. Then we add the $i$th row to any row that has 
$-1$ as an entry below the diagonal 
in the $i$th column. This yields precisely the original matrix $A$, and finishes the proof of the 
first part of Theorem \ref{theo1}.

\smallskip

For the second part of Theorem \ref{theo1} let $A$ be a Wirth matrix 
defining a reflexive cs-crosspolytope $P = \conv(\pm e_1, \ldots, \pm e_d)$. 
Let $i \in \{1, \ldots, d\}$. Now it is 
easy to see that $P$ splits into $R_i := \conv(\pm e_j \,:\, j \not=i)$ and $\conv(\pm e_i)$ if 
and only if the so-called normalized volume of $R_i$ is the same as the normalized volume of $P$. Equivalently, 
the greatest common divisor of all $d-1$-minors of the matrix formed by the columns 
$e_j \,:\, j \not= i$ equals the determinant of $A$. However this happens if and only if the 
$i$th row of $A$ contains only one $1$ (namely on the diagonal):

To see this assume $A_{i,i} = 1$ and $A_{i,j} = 1$ for $i \not= j$. Then we remove the $i$th column and the $j$th 
row. By Laplace developing one gets that this $d-1$-minor 
can never be $\det(A)=2^f$, where $f$ is the number of $2$'s in $A$.

This proves that $A$ is a 1-minimal Wirth matrix if and only if $P$ is 1-irreducible. 
Moreover the splitting of $P$ into a 1-irreducible $P'$ and copies of $[-1,1]$ is unique 
up to isomorphism as lattice polytopes, since the vertices of $P'$ are uniquely determined.

\subsection*{Proof of Theorem \ref{generalewald}}

Let $P$ be a simplicial reflexive polytope with facets $F, -F$. Proposition \ref{simplizchar} immediately yields 
that $P$ splits in a 1-irreducible reflexive cs-crosspolytope $P'$, 
$r$ copies of $[-1,1]$, del Pezzo polytopes $V_{k_1}, \ldots, V_{k_s}$, 
and pseudo-del Pezzo polytopes $\Vs_{p_1}, \ldots, \Vs_{p_t}$, where the 1-irreducible reflexive cs-crosspolytope 
$P'$ is given by the reduction of the Wirth matrix $A$. Hence we have the splitting in Theorem \ref{generalewald}. 

Abbreviating, we let $P$ split into $Q$ and $S$, where $Q$ splits into $P'$ and $r$ copies of $[-1,1]$, and $S$ splits into the remaining factors. 
We define two vertices of $P$ to be connected, if they are contained in a minimal non-face of $P$ 
(a {\em primitive collection} in the language of Batyrev \cite{Bat91}).  From the description of $V_d$ and $\Vs_d$ in \cite{Cas03} 
we see that the connected components of $\V(P)$ of size $> 2$ are precisely the vertex sets of the factors in the splitting of $S$. 
Hence the numbers $k_1, \ldots, k_s, p_1, \ldots, p_t$ are (even combinatorial) invariants of $P$. 
Moreover the isomorphism type of $P$ as a lattice polytope determines the isomorphism type of $Q$, and hence 
$P'$ and $r$ by the second part of Theorem \ref{theo1}.

\section{Applications}

\subsection*{All facets are isomorphic}

\begin{corollary}
Let $P$ be a pseudo-symmetric simplicial reflexive polytope. 

Then any two facets of $P$ are isomorphic as lattice polytopes. 
Especially any two facets of $P$ have the same number of lattice points.
\label{facetsiso}
\end{corollary}

\begin{proof}

Since all facets of a smooth Fano polytope are isomorphic, by Theorem \ref{generalewald} we only have 
to regard reflexive cs-crosspolytopes. Now we use the first part of Theorem \ref{theo1}.

\end{proof}

\subsection*{The maximal number of vertices}

\begin{corollary}
Let $P \subseteq \MR$ be a pseudo-symmetric simplicial reflexive polytope. 

If $d$ is even, then 
$$\abs{\V(P)} \leq 3 d,$$
with equality only if $P$ splits into $d/2$ copies of $V_2$.

If $d$ is odd, then
$$\abs{\V(P)} \leq 3 d - 1,$$
with equality only if $P$ splits into $[-1,1]$ and $(d-1)/2$ copies of $V_2$.

\label{csymmy}
\end{corollary}

\begin{proof}

By Theorem \ref{generalewald} $P$ splits into 
an $l$-dimensional reflexive cs-crosspoly\-tope $Q$ and a smooth Fano polytope $R$ splitting 
into del Pezzo polytopes $V_{k_1}, \ldots, V_{k_s}$ 
and pseudo-del Pezzo polytopes $\Vs_{p_1}, \ldots, \Vs_{p_t}$, so $l + k_1 + \cdots + k_s + p_1 + \cdots + p_t = d$, 
where $k_1, \ldots, k_s, p_1, \ldots, p_t$ are even. 

Hence $\abs{\V(P)} = 2 d + 2 s + t \leq 2 d + \dim(R) -t \leq 3 d$. 

If $\abs{\V(P)} = 3 d$, then $\dim(R) = d$, $P \cong R$, $d$ is even, $t=0$ and $s=d/2$.

So let $d$ be odd and $\abs{\V(P)} = 3 d - 1$. 
Since $P$ has dimension $d$, so is not even-dimensional, we cannot have $\dim(R) = d$, hence 
$\dim(R) = d-1$. Therefore $Q$ is one-dimensional, so isomorphic to $[-1,1]$. Now use the first statement 
for $R$ with $\abs{\V(R)} = 3 (d-1)$.

\end{proof}

In \cite{Cas04} the first part was recently shown to hold for arbitrary simplicial reflexive polytopes. 
The second part is an extension of \cite[Thm. 5.9]{Nil05}.

\subsection*{The maximal number of facets}

\begin{remark}{\rm 
First we look at the (pseudo) del Pezzo polytopes: For this we use the notation of Definition \ref{delpezzos}. 
By setting $e_0 := - e_1 - \cdots - e_d$ it is a straightforward calculation (see 
\cite{VK85} or \cite{Cas03}) that the facets of $V_d$ have as vertices precisely $\{\pm e_j \;:\; 
j = 0, \ldots, d,\, j \not= i\}$, for fixed $i \in \{0, \ldots, d\}$, where 
exactly half of the signs are equal to $+1$ and the others are equal to $-1$. Hence we get:
\[\abs{\F(V_d)} = (d+1) \binom{d}{\frac{d}{2}}.\]
In just the same way (see \cite{Cas03}) we can calculate
\[\abs{\F(\Vs_d)} = d \binom{d-1}{\frac{d}{2}} + \sum\limits_{i=\frac{d}{2}}^{d} \binom{d}{i}.\]
\label{pezzoremark}}
\end{remark}

\begin{corollary}
Let $P \subseteq \MR$ be a pseudo-symmetric simplicial reflexive polytope. 
Then 
\[\abs{\F(P)} \leq 6^{d/2},\]
where equality is only attained if $d$ is even and $P$ splits into $d/2$ copies of $V_2$.
\end{corollary}

\begin{proof}

By Theorem \ref{generalewald} we can assume that $P$ splits into 
an $l$-dimensional reflexive cs-crosspolytope $Q$, del Pezzo polytopes $V_{k_1}, \ldots, V_{k_s}$, 
and pseudo-del Pezzo polytopes $\Vs_{p_1}, \ldots, \Vs_{p_t}$, so $l + k_1 + \cdots + k_s + p_1 + \cdots + p_t = d$, 
where $k_1, \ldots, k_s, p_1, \ldots, p_t$ are even.

Now it is not difficult to show that
\[2n \binom{2n-1}{n} + \sum\limits_{i=n}^{2n} \binom{2n}{i} < (2n+1) \binom{2n}{n} \leq 6^n\]
for $n \in \N_{\geq 1}$, with equality at the right only for $n=1$. 

Hence 
the previous remark and $\abs{\F(Q)} = 2^l \leq 6^{l/2}$ yields 
$\abs{\F(P)} \leq 6^{l/2} 6^{k_1/2}$\\$\cdots 6^{k_s/2} 6^{p_1/2} \cdots 6^{p_t/2}= 6^{d/2}$, 
where equality implies that $d$ is even and $P$ splits into $d/2$ copies of $V_2$.

\end{proof}

Since $V_2^* \cong V_2$, this yields for all centrally symmetric simple reflexive polytopes 
a confirmation of 
the general conjecture in \cite{Nil05} that $V_2^{d/2}$ is the single reflexive polytope with 
the maximal number of vertices $6^{d/2}$.

\subsection*{The maximal number of lattice points}

In \cite[Thm. 6.1]{Nil06} it was shown that $[-1,1]^d$ solely contains the most 
lattice points among all $d$-dimensional centrally symmetric reflexive polytopes. 
Here we prove an analogous result for pseudo-symmetric simplicial reflexive polytopes.

\begin{definition}
Let $D_d$ be the 1-irreducible reflexive cs-crosspolytope associated to 
the 1-minimal Wirth matrix
\[A_d := \begin{pmatrix}2 \id_{d-1} & 0\\ 1 \cdots 1 & 1\end{pmatrix}\]
in a lattice basis $m_1, \ldots, m_d$.
\end{definition}

\begin{remark} 
The dual polytope $D_d^* = \conv(\pm (m^*_d - x) \;:\; x \in \sum_{i=1}^{d-1} c_i m^*_i, \; c_i \in \{0,1\})$, 
in the dual lattice basis of $m_1, \ldots, m_d$, is a 
reflexive polytope where the vertices are the only lattice points on the boundary.
\label{dd}
\end{remark}

\begin{corollary}
Let $P \subseteq \MR$ be a pseudo-symmetric 
simplicial reflexive polytope. Then we have:

\begin{itemize}
\item Any lattice point on $\randp$ is either a vertex or 
the middle point of an edge.
\item $\abs{P \cap M} \leq 2 d^2 + 1$. Any facet of $P$ has at most $\binom{d+1}{2}$ lattice points.
\item The following statements are equivalent:
\begin{enumerate}
\item $\abs{P \cap M} = 2 d^2 + 1$
\item Some facet has $\binom{d+1}{2}$ lattice points
\item Any facet has $\binom{d+1}{2}$ lattice points
\item $P \cong D_d$
\end{enumerate}
\end{itemize}
\label{last}
\end{corollary}

\begin{proof}

Let $F,-F \in \F(P)$. Applying Proposition \ref{simplizchar} to $P$ and $F$ we can assume that $\V(F) = \{e_1, \ldots, e_d\}$ are the 
columns of a Wirth matrix $A$.

To the first point: Since by Corollary \ref{facetsiso} any two facets are isomorphic, we may assume $m \in F \cap M$, 
$m \not\in \V(F)$. By looking at $A$ we get $m = (e_i + e_j)/2$ for a pair $i < j$ with $(e_i)_i = 2$. 
Proceeding further this shows that
\[\abs{F \cap M} \leq \binom{d}{2} + d = \frac{(d+1) d}{2},\]
where equality necessarily implies $A = A_d$, and hence $P = D_d$, since by the fourth point of Proposition \ref{simplizchar} 
we have $\abs{\V(P)} = 2 d$. 
By Corollary \ref{facetsiso} this proves the bound on the number of lattice points in any facet, and $(2) \ra (3) \ra (4)$.

Now we use the notation in Lemma \ref{fund2} and Proposition \ref{simplizchar}, so $u := \eta_F$. 
Let $m \in \randp \cap M \cap u^\perp$. If $m$ is a vertex, then $m \in \{v^1, \ldots, v^d\}$. 
So let  $m \not\in \V(P)$. By Lemma \ref{adjacent} we still have $m \in F_i$ 
for some $i \in \{1, \ldots, d\}$. In particular $v^i \not\in u^\perp$, so $v^i = - e_i$. By the first point of this corollary 
we get $m = (-e_i + e_j)/2$ for some $j \not= i$. 

This yields\\
\centerline{$\abs{\randp \cap M \cap u^\perp} \leq \sum\limits_{i=1}^d \abs{F_i \cap M \cap u^\perp} \leq d (d-1)$.}
Hence, since $P$ contains no non-zero interior lattice points (e.g., \cite[Prop. 1.12]{Nil05}), we have the upper bound
\[\abs{P \cap M} \leq 1 + 2 \abs{F \cap M} + d (d-1) \leq 2 d^2 + 1,\]
where equality implies $\abs{F \cap M}$ to be maximal. This proves the remaining implications $(4) \ra (1) \ra (2)$.

\end{proof}

The proof also shows how to easily read off the number of lattice points of a reflexive cs-crosspolytope from the associated Wirth matrix.

\subsection*{Embedding in the standard lattice cube}

In \cite[Conjecture 2]{Ewa88} Ewald conjectures that{\em, up to unimodular transformation, all vertices of a 
smooth Fano polytope have coordinates $1,-1,0$ only.} 
We say that there is an {\em embedding} into $[-1,1]^d$. This is proven for $d\leq 4$ by the existing classification. 
In general the conjecture does not hold for simplicial reflexive polytopes, even in dimension two. 
However it is true, if we assume pseudo-symmetry:

\begin{corollary}(Wirth, N.)
Let $P$ be a pseudo-symmetric simplicial reflexive polytope. 
Then $P$ can be embedded into $[-1,1]^d$.
\label{embedyes}
\end{corollary}

\begin{proof}

Since (pseudo) del Pezzo polytopes are by definition contained in $[-1,1]^d$, 
by Theorem \ref{generalewald} and \ref{theo1} we just have to show that performing row operations on a Wirth matrix $A$ we 
get a matrix containing only $\{-1,0,1\}$. 
If we assume $A_{j,j} = 2$, then there is an $i > j$ (minimal) such 
that $A_{i,j} = 1$. Then we just have to subtract the $i$th row from the $j$th. 
We proceed by induction on $j$.

\end{proof}

This result in the case of a reflexive cs-crosspolytope can be found in 
\cite[Satz 4.4]{Wir97}. In \cite[Kapitel 4]{Wir97} there is a discussion of 
the topic of embedding, where we can also find the following example that shows that we cannot drop 
the assumption of simpliciality:

\begin{remark}
Let $P \subseteq \MR$ be a centrally symmetric reflexive polytope. Then $P^*$ can be embedded into 
$[-1,1]^d$ 
if and only if $P$ contains a lattice basis of $M$. 

This is, for instance, not true for the four-dimensional reflexive cs-cross\-poly\-tope $P$ 
associated to the 1-minimal Wirth matrix 
$\left(\begin{smallmatrix} 2 & 0 & 0 & 0\\ 1 & 1 & 0 & 0\\ 1 & 0 & 1 & 0\\ 1 & 0 & 0 & 1\end{smallmatrix}\right)$, since 
any lattice point on the boundary is $\pm$ a column.
\end{remark}

\smallskip

However we can still embed $P^*$ into a small multiple of $[-1,1]^d$:

\begin{corollary}
Let $P$ be a pseudo-symmetric simplicial reflexive polytope. 
Then $P^*$ can be embedded into $\lfloor\frac{d}{2}\rfloor [-1,1]^d$.
\label{embeddual}
\end{corollary}

\begin{proof}

Since the duals of the (pseudo) del Pezzo polytopes are always contained in $[-1,1]^d$ (for this use 
Remark \ref{pezzoremark} and \cite{Cas03}), 
by Theorem \ref{generalewald} we can assume $P = \conv(\pm e_1, \ldots, \pm e_d)$. 
Let $m_1, \ldots, m_d$ be the lattice basis of $M$ in Proposition \ref{simplizchar} such that 
$A = \begin{pmatrix}2 \id_{f} & 0\\ C & \id_{d-f}\end{pmatrix} \in \Mat_d(\N)$. Then 
$A^{-1} = \begin{pmatrix}1/2 \id_{f} & 0\\ -C/2 & \id_{d-f}\end{pmatrix}$. Now the rows of $A^{-1}$ are precisely 
the coordinates of the dual $\R$-basis $e_1^*, \ldots, e_d^*$ (in the dual lattice basis $m_1^*, \ldots, m_d^*$). 
Furthermore for any facet $F \in \F(P)$ we have $\eta_F = \pm e_1^* \pm \cdots \pm e_d^* \in N$ 
for some signs $\pm$. Hence the vertices of $P^*$ have coordinates in $[\lfloor -d/2 \rfloor, \lfloor d/2 \rfloor]$ 
with respect to the lattice basis $m_1^*, \ldots, m_d^*$.

\end{proof}

\begin{acknowledgments}
The author would like to thank G. Ewald for giving \cite{Wir97} to his disposal, 
C. Haase for discussing the proof of Theorem \ref{theo1}, and C. Casagrande for pointing out a mistaken statement in 
the introduction. Most results of this paper are part of 
the thesis of the author, advised by V. Batyrev.
\end{acknowledgments}

\bibliographystyle{amsalpha}

\end{document}